\documentclass[12pt, a4paper]{article}

\usepackage{amsfonts, amssymb}
\def\ra{\rightarrow}

\def\Z{\mathbb Z}
\def\Oc{\mathcal O}
\def\R{\mathbb R}

\def\Fc{\mathcal F}

\newtheorem{them}{Theorem}
\newtheorem{lemm}{Lemma}
\title{A  restriction theorem for torsion-free sheaves on some elliptic manifolds}

\author{Victor VULETESCU}

\date{}

\begin{document}
\maketitle

\begin{abstract}
	We prove that if $X$ is the total space of an elliptic principal bundle $\pi:X\ra B$ which is non-k\"ahler, then the restriction of any torsion-free sheaf on $X$ to the general fiber of $\pi$ is semi-stable.
\end{abstract}

\section{Introduction}
In the study of holomorphic vector bundles over a given compact complex manifold $X$, especially in the study of (semi)stable ones,  a very useful tool is the study of their restrictions to general memebers of a given family of subvarieties of $X.$ However, the restriction of a (semi)stable vector bundle to a submanifold is not always semistable. Still, under some strong hypothesis, such as $X$ is projective and the family of subvarieties is a family of divisors "ample enough", the restriction of a stable vector bundle to the general memeber remains (semi)stable: this is "Flenner's restriction theorem", see \cite{Fle}. Flenner's theorem  has been extended to the more general context of algebraic varieties in arbitrary characteristic (see e.g \cite{Langer}), but, to the author's knowledge, there is no such extension to the case of non-projective manifolds. The present note tackles this case.

\section{Notations and basic facts}
The context we are working is the following. We fix a compact complex manifold $B$ and an elliptic curve $F$. To every principal elliptic bundle $$\pi:X\ra B$$ one can associate (up to the obvious action of $SL(2, \Z) $) a couple of elements 
$$\left(c'_1(\pi), c"_1(\pi)\right)\in H^2(B, \Z)\times H^2(B, \Z)$$
called {\em the Chern classes} of the bundle $\pi$ (see e.g. \cite{Br1}).

If at least one of the Chern classes is non vanishing in $H^2(B, \R)$, one can prove by a standard argument using the Leray spectral sequence of the fibration that the homology class of any fiber
$[F]\in H_2(X, \R)$ vanishes; as the fibers are compact complex submanifolds, this shows that $X$ is not of K\"ahler type.

We also recollect the notion of stability; since we will use this concept for vector bundles on curves, we will only recall the definition in this case. Hence, a vector bundle $E$ on a smoothe projective curve  will be called {\em stable} (respectively {\em semistable}) if for any subbundle $\Fc\subset E$ with $0<rank(\Fc)<rank(E)$ one has
$$\frac{deg(\Fc)}{rank(\Fc)}<\frac{deg(E)}{rank(E)}$$
(resp $"\leq"$ for semistability). A vector bundle which is not semistable is called {\em unstable}.

Eventually, let us recall a concept which is of relevance only on non-algebraic complex manifolds. If $X$ is compact complex manifold and $\Fc$ is a coherent sheaf on $X$, then  $\Fc$ is called {\em reducible} if there exist a coherent subsheaf $\Fc' \subset \Fc$ with $0<rank(\Fc')< rank (\Fc);$ if no such subsheaf exist then $\Fc$ is called {\em irreducible.} Notice that on projective manifolds all coherent sheaves are reducible; still, on general compact complex manifolds this is not always the case, as one can see for instance looking at the tangent bundle of a  K3 surface  $X$ with $Pic(X)=0$ (the general K3 surface is so).

%%%%%%%%%%%%%%%%%%%%%%%%%
\section{Some Lemmas}

In the following we  some lemmas, which are most likely classical and well-known; but since we don't have any precise reference, we include the proofs here.

\begin{lemm}\label{lema1}
Let $\pi:X\ra B$ be an elliptic principal bundle. If the homology class $[F]\in H_2(X, \R)$ vanishes (i.e the Poincar\'e dual $PD_F$ is zero), then any proper closed analytic subset $Y\subset X, dim(Y)<dim(X),$ does not meet the general fiber.
\end{lemm}

{\bf Proof.} The only non-obvious case is when $Y$ is a hypersurface. But in this case, if $Y$ meets all the fibers, then it meets the general fiber transversely in finitely many points. But then
$$0<\#(Y\cap F)=\int_XPD_Y\wedge PD_F=0$$
since $PD_F=0$ by the assumption that $0=[F]\in H_2(X, \R).$

\begin{lemm}\label{lema2} For $X$ as in the previous Lemma and 
for any torsion-free sheaf $E$ on $X$ we have $$deg(E_{\vert F})=0$$ for $F=$general fiber of $\pi.$

\end{lemm}

{\bf Proof}. Indeed, as $E$ is torsion-free, we see $Sing(E)$ has codimension at least two. Let $L=det(E)^{\vee \vee}$ be the bidual of the determinant of $E$; it is a reflexive sheaf of rank one on $X$, so it is a line bundle (cf e.g. \cite{OSS}). Moreover, the map $det(E)\ra L$ is an isomorphism outside $sing(E)$, so if $F$ is any fiber not meeting $Sing(E)$ we have

$$deg(E_{\vert F})=deg(det(E)_{\vert F})=deg(L_{\vert F})=i^*(c_1(L))$$
where $i:F\ra X$ is the inclusion of the fiber $F$.
But as $[F]=0$ in $H_2(X, \R)$ we see $i^*(c_1(L))=0,$ Q.E.D. Lemma.

\begin{lemm}\label{lema3}
If $F$ is an elliptic curve  and if $E$ is a vector bundle of degree zero on $F$ which is generated by its global sections, then $E$ is trivial.
\end{lemm}

We use the following argument from L. Ein (cf \cite{Ein}, Proposition 1.1):

\noindent "{\bf Lemma}. If $X$ is a compact complex manifold, and $E$ is a globally generated vector bundle on $E$ such that its dual $E^\vee$ has a section, then $E$ splits as $E=\Oc_X\oplus F.$ "

We do induction of $rank(E).$ For $rank(E)=1$ the assertion is immediate. If $rank(E)\geq 2$, letting $K=Ker\left(H^0(F, E)\otimes {\Oc}_F\ra E\right)$ we get an extension:
\begin{eqnarray}\label{exte}
0
\ra
K
\ra 
H^0(F, E)\otimes {\Oc}_F
\ra 
E
\ra 
0.
\end{eqnarray}
Now, either the extension splits (and hence $E$ is trivial), or 
$$H^1(F, E^\vee\otimes K)\not=0.$$
As $deg(E)=0$ we have also $deg(K)=0$ so we further get by Riemann-Roch on F that
\begin{eqnarray}\label{h0}
H^0(F, E^\vee\otimes K)\not=0.
\end{eqnarray}
Twisting the above extension (\ref{exte}) by $E^\vee$ we get

$$0
\ra
K\otimes E^\vee
\ra 
H^0(F, E)\otimes E^\vee
\ra 
E\otimes E^\vee
\ra 
0$$
hence, from (\ref{h0}), we get 
$$H^0(F, E^\vee)\not=0$$
Applying Ein's Lemma, we get $E=\Oc_F\oplus E_1$. But  $E_1$ has degree zero and is generated by its global sections too, so by the induction hypothesis, $E_1$ is trivial. Consequently,  $E$ is trivial too.

\begin{lemm}\label{sublema}
Let $F$ be an elliptic curve and $L$ a semistable vector bundle on $F$ such that $deg(L)=0$. Then  there is a Zariski-open subset $U\subset Pic_0(F)$ such that
$H^0(F, L\otimes I)=0$ for all $I\in U.$ 
\end{lemm}

{\bf Proof.} (See also \cite{Ray}).
Again, we do induction on $rank(L).$ For $rank(L)=1$ the claim is immediate (take $U=Pic_0(F)\setminus\{L^\vee)\}$), so assume $rank(L)>0$.

In the case $H^0(F, L)=0$, from the existence of the Poincar\'e bundle and Grauert's upper continuity theorem we get $H^0(F, L\otimes I)=0$ for all $I$ in a Zariski neiborghood of $\Oc_F.$

In the case $h^0(F, L)>0$ take some $s\in H^0(F, L), s\not=0;$ it defines a map
$$0\ra \Oc_F\stackrel{s}{\ra}L$$
We infer that this map has torsion-free cokernel; since otherwise, moding out by the torsion of the cokernel, we would get a nontrivial map into $L$ from a nontrivial, effective  divisor on $F$, contradicting the hypothesis that $L$ is semistable. So $L$ sits in an exact sequence

$$
0
\ra
\Oc_F
\ra
L
\ra L'
\ra 
0
$$
with $L'=$torsion-free (hence locally free, as $F$ is a curve); in particular, $deg(L')=0.$
It is easy to see that $L'$ is semistable too, so by the induction hypothesis $H^0(F, L'\otimes I)=0$ for all $I$ is some open subset $U\subset Pic_0(F).$ So
$$H^0(F, L\otimes I)=0$$ for all $I
\in U\setminus \{\Oc_F\},$ Q.E.D. Lemma.

Eventually, we recollect a fact which is true more generally
\begin{lemm}\label{lema4}
Let $F$ be  an elliptic curve and 
$$0
\ra
L
\ra M
\ra R
\ra 0
$$
an exact sequence of vector bundles of $F$ with $$deg(L)=deg(R)=0.$$
If $L$ and $R$ are semistable, them $M$ is semistable too.
\end{lemm}

{\bf Proof.}  Using Lemma \ref{sublema} we get a line bundle $I\in Pic_0(F)$ such that 
$$H^0(F, R\otimes I)=H^0(F, L\otimes I)=0;$$
this implies $H^0(F, M\otimes I)=0$ as well.

So, replacing $M$ by $M\otimes I$ we can further assume
$H^0(F, M)=0$. Now, if $M$ would be unstable, we would get a destabilizing vector subbundle $D\subset M$ with $deg(D)>0$. But $deg(D)>0$ implies $H^0(F, D)\not =0;$ so $H^0(F, M)\not=0$ as well, contradiction, Q.E.D. Lemma.

%%%%%%%%%%

\section{The main result}

We are now in position to state and prove the main result.

\begin{them}
Let $\pi:X\ra B$ be an elliptic principal bundle with at least one of the Chern classes non-vanishing in $H^2(B, \R)$ (in particular, $X$ is nonK\"ahler).
Then the restriction of any torsion-free sheaf $E$ on $X$ to the general fiber of $\pi$ is semi-stable.
\end{them}

Before proving it, let us make a small comment. As one can see, the theorem gives the semi-stability of the restriction of $E$ to the general fiber of $\pi$ with  {\em no apriori assumptions like (semi)stability for} $E.$ This is not completely surprising; in the non-projective context, more exactly on non-projective surfaces, the "Bogomolov inequality" $\Delta(E)\geq 0$, holds similarly for {\em all torsion-free sheaves} $E$ (cf \cite{BaLp}, or \cite{Br2} for a simpler proof), in contrast to the projective case, when it holds mainly for stable vector bundles.

\medskip

{\bf Proof of the theorem.}  We do induction on the rank $r=rk(E).$ For $r=1$ there is nothing to prove, so we assume $r\geq 2.$

\medskip

{\bf Case 1: $E$ is reducible.} That is, $E$ sits in an exact sequence
$$0
\ra
L
\ra
E
\ra
R
\ra 
0$$

By the Lemma \ref{lema2}, we see that for a general fiber $F$ of $\pi,$ $L_{\vert F}, R_{\vert F}$ are locally free of degree zero. More, by the induction hypothesis, both $L_{\vert F}, R_{\vert F}$ are also semistable, so $E_{\vert F}$ is semistable too, by Lemma \ref{lema4}.

{\bf Case 2: $E$ is irreducible}. We distinguish again two subcases:

{\bf Subcase 2.1: $\pi_*(E)=0$.} In this case, $H^0(F, E_{\vert F})=0$ for the general fiber. But as also $deg(E_{\vert F})=0$ for the general fiber $F$, we see at once that $E_{\vert F}$ is semistable. Indeed, if this is not the case, then a destabilizing subsheaf $D\subset E$ would have $deg(D)>0;$ but then $h^0(F, D)>0$ so $h^0(F, E_{\vert F})>0$ too,  contradiction.

{\bf Subcase 2.2: $\pi_*(E)\not=0$.} Let $\alpha:\pi^*\pi_*(E)\ra E$ be the canonical morphism and let $\Fc=Im(\alpha).$ As $E$ is irreducible and as $\alpha$ is non-trivial, we see we have 
$$rank(\Fc)=rank(E).$$
Let $Y=Supp(E/\Fc);$ by Lemma \ref{lema1}, $Y$ cannot meet all the fibers of $\pi$ so for the general fiber $F$ we have $\Fc_{\vert F}=E_{\vert F};$ more, by Lemma \ref{lema2} we can assume $deg(E_{\vert F})=0.$

So, for the general fiber $F$ we have a surjection
$$\pi^*\pi_*(E)_{\vert F}\ra E_{\vert F}.$$

\noindent But 
$$\pi^*\pi_*(E)_{\vert F}$$ is trivial, so $E_{\vert F}$ is spanned by its global sections. As it is also of degree zero, it folllows by Lemma \ref{lema3} that $E_{\vert F}$ is trivial, in particular semi-stable.

\medskip
{\bf Acknowledgements.} The author would like to thank V. Br\^{\i}nz\u anescu for asking the question of the semistability of restrictions in the non-K\"ahler context and for a careful reading of some firsts drafts of the paper. The ellaboration of paper was funded by the grant "Vector Bundle Techniques in the Geometry of Complex Varieties", PN-II-ID-PCE-2011-3-0288,  Contract 132/05.10.2011.

\noindent {\small UNIVERSITATEA BUCURE\c STI, FACULTATEA DE MATEMATIC\u A \c SI INFORMATIC\u A, and\\
"SIMION STOILOW" INSTITUTE OF MATHEMATICS OF THE ROMANIAN ACADEMY\\
%email: vuli\@fmi.unibuc.ro}
%Universitatea Bucure\c sti, Facultatea de Matematic\u a \c si Informatic\u a, Romania, and\\
%"Simion Stoilow" Institute of Mathematics of the Romanian Academy, Romania }


\begin{thebibliography}{100}

\bibitem{BaLp}
 C. B\u anic\u a; J.  Le Potier, {\em Sur l'existence des fibr\'es vectoriels holomorphes sur les surfaces non- alg\'ebriques. (On the existence of holomorphic vector bundles on non- algebraic surfaces). }
J. Reine Angew. Math. 378, 1-31 (1987).

\bibitem{Br1} V. Br\^{\i}nz\u anescu,
{\em Neron-Severi group for nonalgebraic elliptic surfaces. I: Elliptic bundle case. }
Manuscr. Math. 79, No.2, 187-195 (1993).

\bibitem{Br2} V. Br\^{\i}nz\u anescu, 
{\em A simple proof of a Bogomolov type inequality in the case of nonalgebraic surfaces. }
Rev. Roum. Math. Pures Appl. 38, No.7-8, 631-633 (1993).

\bibitem{Ein} L. Ein, {\em An analogue of Max Noether's theorem. }
Duke Math. J. 52, 689-706 (1985).

\bibitem{Fle}H. Flenner, {\em  Restrictions of semistable bundles on projective varieties,} Comment. Math. Helv. 59 (1984), 635–650.

\bibitem{Langer} A. Langer, {\em A note on restriction theorems for semistable sheaves}, Math.Res.Lett.17(2010),no.05,823–832

\bibitem{OSS}  Ch. Okonek; M. Schneider; H. Spindler,
{\em Vector bundles on complex projective spaces. }
Progress in Mathematics. 3. Boston - Basel - Stuttgart: Birkhäuser. VII, 389 p.

\bibitem{Ray} M. Raynaud, {\em Sections des fibr\'es vectoriels sur une courbe}, Bulletin de la S.M.F., tome 110 (1982), p.103-125

\end{thebibliography}
\end{document}